\documentclass[11pt,oneside]{article}
\usepackage{amssymb,amsmath}
\usepackage[english]{babel}
\usepackage{url}
\renewcommand{\le}{\leqslant}
\renewcommand{\ge}{\geqslant}

\newcommand{\rk}{\mathbf{rk}}

\newcommand{\RR}{\mathbb{R}}
\newcommand{\ZZ}{\mathbb{Z}}

\newcommand{\NN}{\mathbb{N}}

\newcommand{\SSS}{\mathcal{S}}

\newcommand{\RRR}{\mathbf{R}}
\newcommand{\vz}{\mathbf{z}}

\newcommand{\vx}{\mathbf{x}}

\newcommand{\vy}{\mathbf{y}}
\newcommand{\va}{\mathbf{a}}
\newcommand{\vb}{\mathbf{b}}

\newcommand{\vs}{\mathbf{s}}
\newcommand{\vc}{\mathbf{c}}

\newtheorem{lemma}{Lemma}
\newtheorem{theorem}{Theorem}
\newtheorem{proposition}{Proposition}

\newtheorem{example}{Example}

\newtheorem{corollary}{Corollary}

\newtheorem{theoremls}{Theorem LS}

\newenvironment{proof}{\textsc{Proof. }}{\ \newline\hspace*{\fill}$\boxtimes$}

\newcommand{\vA}{\mathbf{A}}
\newcommand{\vB}{\mathbf{B}}

\newcommand{\vM}{\mathbf{M}}

\begin{document}

\title{Upper bounds for the uniform simultaneous Diophantine exponents}

\author{
 Dzmitry Badziahin
}

\maketitle

\begin{abstract}
We give several upper bounds for the uniform simultaneous
Diophantine exponent $\widehat{\lambda}_n(\xi)$ of a transcendental
number $\xi\in\RR$. The most important one relates
$\widehat{\lambda}_n(\xi)$ and the ordinary simultaneous exponent
$\omega_k(\xi)$ in the case when $k$ is sufficiently smaller than
$n$. In particular, in the generic case $\omega_k(\xi)=k$ with a
properly chosen $k$, the upper bound for $\widehat{\lambda}_n(\xi)$
becomes as small as $\frac{3}{2n} + O(n^{-2})$ which is
substantially better than the best currently known unconditional
bound of $\frac{2}{n} + O(n^{-2})$. We also improve an unconditional
upper bound on $\widehat{\lambda}_n(\xi)$ for even values of $n$.
\end{abstract}

{\footnotesize{

Math Subject Classification 2010: 11J13, 11J82}}

\section{Introduction}

Given a real number $\xi$, its {\it $n'$th (uniform) simultaneous
Diophantine exponent} $\lambda_n(\xi)$ (respectively,
$\widehat{\lambda}_n(\xi)$) is defined as the supremum of all
$\lambda$ such that the system of inequalities
$$
\max_{1\le i\le n} |x_0\xi^i - x_i|<Q^{-\lambda},\quad 1\le |x_0|\le
Q;
$$
has a solution $\vx = (x_0,x_1,\ldots, x_n)\in\ZZ^{n+1}$ for an
increasing unbounded sequence of values $Q$ (respectively, for all
large enough values of $Q$).

The {\it $n'$th (uniform) dual Diophantine exponent} $\omega_n(\xi)$
(resp., $\widehat{\omega}_n(\xi)$) is defined as the supremum of
$\omega$ such that the system of inequalities
$$
0<|P(\xi)|\le Q^{-\omega},\quad \deg(P)\le n,\quad H(P)\le Q
$$
has infinitely many solutions $P\in\ZZ[x]$ for an unbounded sequence
of values $Q$ (respectively, for all large values $Q$). Here, $H(P)$
is the trivial height of the polynomial $P$, i.e. it is the supremum
norm of the vector of its coefficients.

The Diophantine exponents are extensively studied in the literature.
They help us understand the approximational properties of real
numbers. The notion of $\omega_n(\xi)$ was first introduced by
Mahler in 1932~\cite{mahler_1932} where he defined his
classification of real numbers. To the best of my knowledge, the
simultaneous Diophantine exponents were formally defined much later
(perhaps, in~\cite{bug_laur_2005}) but were studied long before that
without having an explicit name.

The classical Dirichlet and Minkowski theorems imply that for
transcendental $\xi$, $\lambda_n(\xi)\ge \widehat{\lambda}_n(\xi)\ge
1/n$ and $\omega_n(\xi)\ge \widehat{\omega}_n(\xi)\ge n$ for all
$\xi\in\RR$. It was conjectured by Mahler and later established by
Sprind\v{z}uk~\cite{sprindzuk_1969} that these bounds are in fact
equalities for almost all values of $\xi$ in terms of Lebesgue
measure. The precise values of all Diophantine exponents are also
known for algebraic numbers, thanks to Schmidt Subspace
Theorem~\cite{schmidt_1980}. If $\xi$ is a real algebraic number of
degree $d$ then
$$
\omega_n(\xi) = \omega_n(\xi) = \min\{n,d-1\};\quad \lambda_n(\xi) =
\widehat{\lambda}_n(\xi) = \frac{1}{\min\{n,d-1\}}.
$$
Further in this paper, we will always assume that $\xi$ is
transcendental.

Diophantine exponents play an important role in understanding, how
well a given real number $\xi$ can be approximated by algebraic
numbers of bounded degree. Heuristically, it is easy to see for
$\omega_n(\xi)$. Indeed, large values of $\omega_n(\xi)$ imply that
there exist polynomials that take very small values at $\xi$, and
that in turn implies that one of the roots of those polynomials is
near $\xi$. This idea was used by Wirsing~\cite{wirsing_1961} in
1961 where he showed that for any transcendental $\xi\in\RR$ and
$\omega^* < \frac{\omega_n(\xi)+3}{2}$ there exist infinitely many
algebraic real numbers $\alpha$ of degree at most $n$ such that
$$
|\xi - \alpha| < H(\alpha)^{-\omega^*}.
$$
Later, this inequality was strengthened in~\cite{tsishchanka_2007,
ber_tsi_1994, bug_sch_2016}. In 1969, Davenport and
Schmidt~\cite{dav_sch_1969} discovered a relation between $\omega^*$
and $\widehat{\lambda}_n(\xi)$. Wirsing himself conjectured that for
all transcendental $\xi$ the degree $\omega^*$ can be made
arbitrarily close to $n+1$. This problem remains open for all
dimensions $n\ge 3$. The best known bounds on $w^*$ can be found
in~\cite{bad_sch_2020}. For more properties and relations between
the Diophantine exponents of $\xi$ and the notion of $\omega^*$, we
refer to an extensive overview~\cite{bugeaud_2014}.

In this paper, we will closely look at the exponents
$\widehat{\lambda}_n$. We do not know much about what values can
they take, as well as the uniform simultaneous exponents
$\widehat{\omega}_n$. It is known that for $n=2$,
$\widehat{\lambda}_2(\xi)$ can change between $\frac12$ and
$\frac{\sqrt{5}-1}{2}$ where both upper and lower bounds are
sharp~\cite{roy_2003}. Recent developments about the spectrum of
values $\widehat{\lambda}_2(\xi)$ and $\widehat{\omega}_n(\xi)$ can
be found in~\cite{poels_2021, bug_laur_2005}. However for $n\ge 3$
and transcendental $\xi$, it is not even known if
$\widehat{\lambda}_n(\xi)$ can take any other value than $1/n$.

In 1969, Davenport and Schmidt~\cite{dav_sch_1969} verified that
$\widehat{\lambda}_n(\xi) \le 1/\lfloor\frac{n}{2}\rfloor$ which, in
view of the lower bound $1/n\le \widehat{\lambda}_n(\xi)$ leaves a
small window for the values of the exponent $\widehat{\lambda}_n$.
Later, their bound was slightly improved by several authors, most
recent results belong to Laurent~\cite{laurent_2003} and
Schleischitz~\cite{schleischitz_2021}. However, all these upper
bounds are still of the form $\frac{2}{n} + O(n^{-2})$.

\begin{theoremls}
For any transcendental $\xi\in\RR$ and integer $n\ge 2$ one has
$$
\widehat{\lambda}_n(\xi)\le \left\{ \begin{array}{ll}
\frac{2}{n+1}&\mbox{for odd $n$}\\
\tau_n&\mbox{for even $n$,}
\end{array}\right.
$$
where $\tau_n$ is the solution in the interval $[\frac{2}{n+2},
\frac{2}{n})$ of the equation
$$
\left(\frac{n}{2}\right)^n x^{n+1} - \left(\frac{n}{2}+1\right)x +
1=0.
$$
\end{theoremls}

In this paper, we improve the upper bound of
$\widehat{\lambda}_n(\xi)$ for even $n$, while asymptotically it is
still in the range $\frac{2}{n} + O(n^{-2})$.

\begin{theorem}\label{th3}
For any positive integer $n$ one has
$$
\widehat{\lambda}_{2n}(\xi)\le t_n = \frac{\sqrt{n^2+4n}-n}{2n}.
$$
where $t_n$ is the positive root of the equation $nx^2 + nx - 1=0$.
\end{theorem}

Observe that for $n=1$ the bound in Theorem~\ref{th3} is the same as
in the result of Davenport and Schmidt~\cite{dav_sch_1969}.
Therefore it is sufficient to prove this theorem for $n\ge 2$. For
larger values of $n$, we get $t_2\approx 0.366 < 0.371\approx
\tau_4$; $t_3\approx 0.264 < 0.268\approx \tau_6$. In general, one
can check that $t_n < \tau_{2n}$ for all $n\ge 2$.

The Diophantine exponent $\widehat{\lambda}_3$ was studied in more
detail by Roy~\cite{roy_2008} where he got a better upper bound
$\widehat{\lambda}_3(\xi) \le \lambda_0\approx 0.4245$ than for
general $n$.

The most important results of this paper show that the upper bounds
on $\widehat{\lambda}_n(\xi)$ can be made much tighter if some of
the Diophantine exponents $\omega_k(\xi)$ with $2k+1\le n$ are close
to their generic values $k$.

\begin{theorem}\label{th1}
Let $\xi\in\RR$ be a transcendental number. Assume that for a given
$k\in\NN$,
\begin{equation}\label{th1_eq1}
\delta_k:=\frac{k}{\omega_k(\xi) + 1-k} \ge 1.
\end{equation}
Then one has
$$
\widehat{\lambda}_n(\xi) \le \left\{ \begin{array}{ll}
\displaystyle \frac{1}{n-k}&\mbox{for}\ 2k+1\le n<
2k+1+\delta_k\\[2ex]
\displaystyle \min\left\{\frac{1}{n -
\left\lceil\frac{n-\delta_k-1}{2}\right\rceil},
\frac{1}{\left\lfloor \frac{n-\delta_k-1}{2}\right\rfloor
+1+\delta_k}\right\}&\mbox{for}\ n\ge 2k+1+\delta_k.
\end{array}\right.
$$
\end{theorem}

While for fixed $k$ and $n$ tending to infinity, the upper bound in
Theorem~\ref{th1} is still of the form $\frac{2}{n}+O(n^{-2})$, it
may become as small as $\frac{3}{2n} + O(n^{-2})$ if $n = 3k+O(1)$
and $\omega_k(\xi) = k$. We manage to slightly improve the bound in
Theorem~\ref{th1} for some values of $n$ and $k$. However, due to
the complicated nature of those bounds, we do not provide them in
the introduction, but leave them in Section~\ref{sec6}, see
Theorem~\ref{th2}.

{\bf Remark.} In fact, the notion $\delta_k$ in Theorem~\ref{th1}
can be replaced by in some cases smaller notion
$\widehat{\omega}_{k,k+1}(\xi)$ which is defined as the supremum of
$\omega$ such that the inequalities
$$
0<|P(\xi)|\le Q^{-\omega},\quad \deg(P)\le k,\quad H(P)\le Q
$$
have $k+1$ linearly independent solutions $P\in\ZZ[x]$ for all large
values of $Q$. While it is not hard to adapt the proof of
Theorem~\ref{th1} to this new notion, we leave the details of the
proof to the interested reader.

The last result of this paper gives a slight improvement of
Theorem~\ref{th1} for the case $n=3$. It provides a better bound on
$\widehat{\lambda}_3(\xi)$ than that of Roy in the case when
$\omega_1(\xi)$ is close to 1.

\begin{theorem}\label{th3}
For a transcendental number $\xi\in\RR$ the Diophantine exponents
$\widehat{\lambda}_3(\xi)$ and $\omega_1(\xi)$ satisfy the
inequality
\begin{equation}\label{th3_eq}
2\widehat{\lambda}_3(\xi)^3 \omega_1(\xi) -
\widehat\lambda_3(\xi)^2(\omega_1(\xi)-1) + 2\widehat\lambda_3(\xi)
\le 1.
\end{equation}
\end{theorem}

For $\omega_1(\xi)=1$ (which is satisfied for almost all $\xi$) this
inequality is equivalent to $\widehat{\lambda}_3(\xi)\le \lambda_3 =
0.42385\ldots$ which is better than the best known bound
$0.4245\ldots$ from~\cite{roy_2008}. In fact,
condition~\eqref{th3_eq} gives a better estimate on
$\widehat{\lambda}_3(\xi)$ than the currently known one for
$\omega_1(\xi) \le 1.07...$.

%

\section{Notation}

For real numbers $A$ and $B$, we write $A\ll B$ if $A\le cB$ for
some absolute constant $c>0$. The notation for $A\gg B$ and $A\asymp
B$ is defined similarly. In the further discussion, we fix a real
number $\xi$ and assume that $\xi\asymp 1$, i.e. $\xi$ is bounded
from above and below by some absolute positive constants.

We will borrow much of the notation from \cite{roy_2008}. Given a
vector $\vx=(x_0,\ldots, x_n)\in\RR^{n+1}$, we define a function
$L(\vx)$ as follows:
\begin{equation}\label{def_l}
L(\vx):=\max_{1\le i\le n} |x_0\xi^i - x_i|.
\end{equation}
This definition implies that for all $k,l\in \ZZ_{\ge 0}$ with
$l+k\le n$ one has
\begin{equation}\label{eq17}
|x_k\xi^l - x_{k+l}| = |(x_k - x_0 \xi^k)\xi^l - (x_{k+l} -
x_0\xi^{k+l})|\ll L(\vx).\
\end{equation}

By $(\vx_i)_{i\in\NN}$ we denote the sequence of minimal points for
$(\xi, \xi^2, \ldots, \xi^n)$. This is a sequence of points $\vx_i =
(x_{i,0}, x_{i,1},\ldots, x_{i,n})$ in $\ZZ^{n+1}$ such that the
following conditions are satisfied:
\begin{itemize}
\item the positive integers $X_i:= ||\vx_i||$ form a strictly
increasing sequence with $X_1=1$;
\item the positive real numbers $L_i:= L(\vx_i)$ form a strictly
decreasing sequence;
\item if some non-zero point $\vx\in \ZZ^{n+1}$ satisfies
$L(\vx)<L_i$ then $||\vx||\ge X_{i+1}$.
\end{itemize}

One can verify that for transcendental $\xi$ the sequence
$(\vx_i)_{i\in\NN}$ is uniquely defined, apart from probably the
first term $\vx_1$. In further discussion we always assume that
$\xi$ is transcendental. Note that since $\xi\asymp 1$ then for
large enough $i$ we have $x_{i,0}\asymp x_{i,1}\asymp\cdots\asymp
x_{i,n}\asymp X_i$.

We fix a positive real number $\lambda$ for which there exists a
constant $c>0$ such that for all $X\ge 1$ the inequalities
$$
||\vx||\le X;\quad L(\vx)\le cX^{-\lambda}
$$
have a solution $\vx\in\ZZ^{n+1}\setminus \{\mathbf{0}\}$. Clearly,
any number smaller than $\hat{\lambda}_n(\xi)$ satisfies this
condition, while $\hat{\lambda}_n(\xi)$ itself may or may not
satisfy it.

By choosing $X=X_{i+1}-\frac12$ we immediately verify that $L_i\le
c(X_{i+1}-\frac12)^{-\lambda}$ or for simplicity we will write it in
a shorter form:
\begin{equation}\label{eq24}
L_i\ll X_{i+1}^{-\lambda}.
\end{equation}


We fix positive real numbers $\omega_k$, $k\in\{1,\ldots, n\}$ for
which the inequalities
\begin{equation}\label{def_omk}
||\vx||\le X;\quad |x_0 + x_1\xi + x_2\xi^2+\cdots + x_k\xi^k| \ll X^{-\omega_k}
\end{equation}
may have solutions $\vx\in\ZZ^{k+1}\setminus\{\mathbf{0}\}$ only for
a bounded set of values $X$. It is straightforward to check that any
number $\omega_k$ bigger than $\omega_k(\xi)$ satisfies this
condition, while $\omega_k(\xi)$ itself may or may not satisfy it.

One can easily verify that any two consecutive vectors $\vx_{i-1}$ and
$\vx_i$ are linearly independent. By $W_i$ we denote the span of these two
vectors:
$$
W_i:= \langle \vx_{i-1}, \vx_i\rangle.
$$

Now we note that if for all $i\ge i_0$ the vectors $\vx_i$ belong to
a proper subspace of $\RR^{n+1}$ then $\xi$ must be an algebraic
number of degree at most $n$. Indeed, in that case without loss of
generality we may assume that the points lie in a hyperplane
$\va\cdot \vx_i = 0$ for some vector $\va$. Since all $\vx_i$ are
integer vectors, $\va$ may be chosen to be integer too. Then the
points $\frac{\vx_i}{x_{i,0}}$ also lie in the same hyperplane.
Finally, we observe that their limit
$$
\lim_{i\to \infty} \frac{\vx_i}{x_{i,0}} = (1,\xi,\xi^2,\ldots,
\xi^n)
$$
is also in that hyperplane and hence $\xi$ is algebraic.

The observation above implies that for $n\ge 2$ there exists
infinitely many indices $i$ such that $\vx_{i-1}, \vx_i$ and
$\vx_{i+1}$ are linearly independent. By $I$ we denote the set of
indices which satisfy this property. Formally,
\begin{equation}
I:= \{i\in\NN: W_i\neq W_{i+1}\}.
\end{equation}

For $i\in I$, by $U_i$ we denote the span of three consecutive
vectors
$$
U_i:= \langle \vx_{i-1}, \vx_i, \vx_{i+1}\rangle = W_i + W_{i+1}.
$$

For $n\ge 3$ and transcendental $\xi$ we must have that the sequence
of subspaces $U_i$ is not eventually constant. In other words, there
must be infinitely many indices $j\in I$ whose successor $i\in I$
satisfies $U_j\neq U_i$. We denote such a set of indices $j$ by $J$.

For larger values of $n$ one can continue defining systems of higher
dimensional subspaces, analogously to $W_i$ and $U_j$, but we will
not use them in this paper.

Given a vector $\vx \in\RR^{n+1}= (x_0, x_1,\ldots, x_n)$, we define
vectors $\vx^{(k,l)}\in\RR^{l+1}$, $k\ge 0$, $k+l\le n$ as follows:
\begin{equation}
\vx^{(k,l)}:=(x_k, x_{k+1}, \ldots, x_{k+l}).
\end{equation}

\section{Preparatory results}

One of the most important ideas introduced by Davenport and
Schmidt~\cite{dav_sch_1969} is to consider vectors of the form
$\vx_i^{(j,m)}$ for a fixed $m<n$ and all $j\in\{0,\ldots, n-m\}$
$i\in \NN$. They are all almost parallel to the vector $(1,
\xi,\xi^2,\ldots, \xi^m)$ therefore the wedge product of such
vectors must have a very small norm.

\begin{proposition}\label{prop2}
Let $m,d\in\ZZ$ satisfy $0\le m\le n$ and $0\le d\le m+1$. Let
$(l_1, k_1), \ldots, (l_d, k_d)$ be a sequence of integer pairs such
that $0<l_1\le l_2\le \cdots \le l_d$ and $0\le k_i \le n-m$ for all
$i\in\{1,\ldots, d\}$. Then
\begin{equation}\label{prop2_eq1}
\left|\left|\bigwedge_{i=1}^d \vx_{l_i}^{(k_i,m)}\right|\right| \ll
X_{l_d} \prod_{i=1}^{d-1} L_{l_i}.
\end{equation}
\end{proposition}

\proof Consider any $d\times d$ minor of the matrix which is
composed of $\vx_{l_1}^{(k_1,m)}, \ldots, \vx_{l_d}^{(k_d,m)}$ as
rows:
$$
M := \left(\begin{array}{cccc} x_{l_1, k_1+j_1}& x_{l_1,
k_1+j_2}&\cdots & x_{l_1, k_1+j_d}\\
x_{l_2, k_2+j_1}& x_{l_2,
k_2+j_2}&\cdots & x_{l_2, k_2+j_d}\\
\vdots&\vdots&\ddots&\vdots\\
x_{l_d, k_d+j_1}& x_{l_d,
k_d+j_2}&\cdots & x_{l_d, k_1+j_d}\\
\end{array}\right)
$$
In order to prove the proposition we need to verify that $|\det M|
\ll X_{l_d} \prod_{i=1}^{d-1} L_{l_i}$.

For each $i$ between 2 and $d$ we multiply the first column of $M$
by $\xi^{j_i}$ and subtract it from $i$'th row of $M$. As a result,
we get a matrix $M^*$ such that, by~\eqref{eq17}, its entries in row
$i$ and columns 2 to $d$ are $\ll L_{l_i}$. Now we expand the
determinant of $M^*$ by the first column to get $|\det M| \ll
X_{l_d} \prod_{i=1}^{d-1} L_{l_i}$. \endproof

The idea of this proposition is that for appropriately chosen
vectors $\vx_{l_i}^{(k_i,m)}$ and for $\lambda$ large enough, the
right hand side of~\eqref{prop2_eq1} becomes smaller than 1. Since
the left hand side is an integer vectors, this implies that
$\vx_{l_i}^{(k_i,m)}$ are linearly dependent. The remaining part of
the arguments is then to verify that linear dependence of those
vectors implies that $\xi$ is algebraic.

The following lemma is proven in~\cite[Lemma~3.1]{bad_sch_2020} and
is based on the ideas of Laurent~\cite{laurent_2003}.

\begin{lemma} \label{lem12}
    Let $m$ be an integer such that $1\leq m\leq n/2$.
    Assume that
    \begin{equation}  \label{eq:only}
    \lambda > \frac{1}{n-m+1}.
    \end{equation}
    Then for any large $i$ the vectors $\vx_i^{(0,n-m)},
    \vx_i^{(1,n-m)},\ldots, \vx_i^{(m,n-m)}$ are linearly
    independent.
\end{lemma}

The next result states that sometimes we can find even more vectors
that are linearly independent. It may be of independent interest.

\begin{proposition}\label{prop1}
Let $m$ and $k$ be nonnegative integers such that either $k\le 2$
and $1\le m+k\le n/2$ or $k\ge 3$ and $1\le m+2k-2\le  n/2 $. Assume
that
$$
\lambda> \left\{\begin{array}{rl} \displaystyle\frac{1}{n-m-k+1}&
\mbox{for}\;
k\le 2\\[2ex]
\displaystyle \frac{1}{n-m-2k+3}& \mbox{for}\; k\ge 3
\end{array}\right.
$$
Let $\vy_0, \vy_2,\ldots, \vy_k$ be large enough linearly
independent minimal points of $\xi$, not necessarily consecutive.
Then the span of $(\vy_i^{(j,n-m)})_{0\le i\le k;\ 0\le j\le m}$ is
at least $m+1+k$.
\end{proposition}

Notice that Lemma~\ref{lem12} is a particular case of this
proposition for $k=0$. We believe that the proposition can be
strengthened so that $\lambda > (n-m-k+1)^{-1}$  should be
sufficient for $k\ge 2$ too but we do not see an easy way of proving
it.

\proof By Lemma~\ref{lem12}, we immediately get that all the vectors
$\vy_0^{(0,n-m)},\ldots, \vy_0^{(m,n-m)}$ are linearly independent.
We will prove a slightly stronger statement than in the proposition:
there exist indices $j_1, j_2,\ldots, j_k$ such that the system of
vectors $(\vy_0^{(i,n-m)})_{0\le i\le m}$, $\vy_1^{(j_1, n-m)},
\ldots, \vy_k^{(j_k,n-m)}$ is linearly independent. We prove this
statement by induction on $k$, where the base $k=0$ is already
verified by Lemma~\ref{lem12}.

If $1\le k\le 2$ we consider the pair $(m+1,k-1)$ in place of $m$
and $k$ together with vectors $\vy_0 ,\ldots, \vy_{k-1}$ and apply
the proposition. By inductional assumption, there exist indices
$j_1,\ldots, j_{k-1}$ such that
$$
(\vy_0^{(i,n-m-1)})_{0\le i\le m+1},\ \vy_1^{(j_1, n-m-1)}, \ldots,
\vy_{k-1}^{(j_{k-1},n-m-1)}
$$
are linearly independent. For $k=2$, there always exist either a
vector $\vy_1^{(j_1-1, n-m)}$ or a vector $\vy_1^{(j_1, n-m)}$. We
assume that the last case is satisfied, the first case can be dealt
analogously. For $k=1$, no cases need to be considered because the
only vectors in the system are of the form $\vy_0^{(i,n-m-1)}$.
Then, by adding one extra component to each of the above vectors,
and throwing away $\vy_0^{(m+1,n-m-1)}$, we get that the system of
vectors
$$
\SSS:= \{(\vy_0^{(i,n-m)})_{0\le i\le m},\ \vy_1^{(j_1, n-m)},
\ldots, \vy_{k-1}^{(j_{k-1},n-m)}\} =\{\vs_1,\ldots, \vs_{k+m}\}
$$
is also linearly independent.

If $k\ge 3$ we consider the pair $(m+1, k-1)$ in place of $m$ and
$k$ together with vectors $\vy_0^{(0,n-1)}, \ldots, \vy_{k-1}^{(0,
n-1)}$. Then the conditions of Proposition~\ref{prop1} are satisfied
for these new parameters and by inductional assumption there exist
indices $0\le j_1, \ldots, j_{k-1}<m$ such that
$$
(\vy_0^{(i,n-m-1)})_{0\le i\le m},\ \vy_1^{(j_1, n-m-1)}, \ldots,
\vy_{k-1}^{(j_{k-1},n-m-1)}
$$
are linearly independent. Then as before, we add one extra component
to each of the above vectors from the right and get that the system
of vectors
$$
\SSS:= \{(\vy_0^{(i,n-m)})_{0\le i\le m},\ \vy_1^{(j_1, n-m)},
\ldots, \vy_{k-1}^{(j_{k-1},n-m)}\} =\{\vs_1,\ldots, \vs_{k+m}\}
$$
is also linearly independent.

Consider the set $\RRR$ of vectors $\vz\in \RR^{n+1}$ such that all
vectors $\vz^{(i, n-m)}$, $i\in \{0,\ldots, m\}$ belong to the span
of $\SSS$. It is straightforwardly verified that this set is a
vector subspace of $\RR^{n+1}$. We will bound its dimension.

For each vector $\vz\in \RRR$ we get
$$
\vz^{(i,n-m)} = \sum_{j=1}^{k+m} a_{i,j} \vs_j.
$$
By the construction of vectors $\vz^{(i,n-m)}$, we get relations for
all $i\in\{0,\ldots, m-1\}, l\in \{1, \ldots, n-m\}$:
\begin{equation}\label{prop1_eq1}
\sum_{j=1}^{k+m} a_{i,j} s_{j, l} = \sum_{j=1}^{k+m} a_{i+1,j}
s_{j,l-1}.
\end{equation}
Compose a vector $\va$ by arranging the coefficients $a_{i,j}$ in
the following way
$$
\begin{array}{rl}
\va:=& ( a_{0,1}, a_{0,2}, \ldots, a_{0,m+1}, a_{1,1},\ldots,
a_{1,m+1},\ldots, a_{m,m+1},\\
& a_{0,m+2}, \ldots, a_{0,m+k}, \ldots, a_{1,m+2}, \ldots,
a_{1,m+k},\ldots a_{m,m+k}).
\end{array}
$$
Then the conditions~\eqref{prop1_eq1} can be written in a matrix
form $\mathbf{M} \va = \mathbf{0}$ where $\vM = (\vA\mid \vB)$. The
part $\vA$ of $\vM$ is written in the block form as
$$
\vA:= \left(
\begin{array}{ccccccc}
H^+&-H^-& \mathbf{0} & \mathbf{0} & \mathbf{0} & \cdots & \mathbf{0}\\
\mathbf{0} & H^+& -H^- & \mathbf{0} & \mathbf{0} & \cdots & \mathbf{0}\\
\mathbf{0} & \mathbf{0} & H^+ & -H^- & \mathbf{0} & \cdots & \mathbf{0}\\
\vdots & \vdots & \vdots & \ddots & \ddots & \ddots & \vdots\\
\mathbf{0} & \mathbf{0} & \mathbf{0} & \cdots & \cdots & H^+ & -H^-
\end{array}
\right)
$$
where
$$
H^+:= ( \vy_0^{(1,n-m-1)}\;\;
\vy_{0}^{(2,n-m-1)}\;\;\ldots\;\;\vy_0^{(m+1,n-m-1)} ),
$$$$
H^-:= ( \vy_0^{(0,n-m-1)}\;\;
\vy_{0}^{(1,n-m-1)}\;\;\ldots\;\;\vy_0^{(m,n-m-1)} )
$$
and vectors $\vy_0^{(i, n-m-1)}$ are considered as columns. The part
$\vB$ is composed similarly
$$
\vB := \left(
\begin{array}{cccccc}
K^+&-K^-& 0 & 0 &  \cdots & 0\\
0 & K^+& -K^- & 0  & \cdots & 0\\
\vdots & \vdots & \ddots & \ddots & \ddots & \vdots\\
0 & 0 & \cdots & \cdots & K^+ & -K^-
\end{array}
\right),
$$
$$
K^+:= ( \vy_1^{(j_1+1,n-m-1)}\;\;
\vy_2^{(j_2+1,n-m-1)}\;\;\ldots\;\;\vy_{k-1}^{(j_{k-1}+1,n-m-1)} ),
$$$$
K^-:= ( \vy_1^{(j_1,n-m-1)}\;\;
\vy_2^{(j_2,n-m-1)}\;\;\ldots\;\;\vy_{k-1}^{(j_{k-1},n-m-1)} ).
$$

Perform the following column operation for $\vA$ which do not change
the rank of this matrix: starting from $i=1$ till $i=m$, subtract
columns $(m+1)(i-1) + i, (m+1)(i-1) + i + 1, \ldots, (m+1)(i-1) + m$
from columns $(m+1)i + i+1, \ldots (m+1)(i+1)$ respectively. Then
going backwards from $i=m$ till $i=2$, subtract columns $(m+1)i + 2,
(m+1)i+3,  \ldots, (m+1)i + i$ from columns $(m+1)(i-1) + 1, \ldots,
m(i-1)$. After that matric $\vA$ transforms to
$$
\vA^* = \left(
\begin{array}{cccccc}
H & \mathbf{0}&\mathbf{0} & \cdots & \mathbf{0} & \mathbf{0}\\
\mathbf{0} & H&\mathbf{0} & \cdots & \mathbf{0} & \mathbf{0}\\
\vdots& \vdots& \ddots& \ddots&\vdots&\vdots\\
\mathbf{0}& \mathbf{0}& \mathbf{0} & \cdots & H & \mathbf{0}
\end{array}
\right)
$$
where
$$
H:= ( \vy_0^{(1,n-m-1)}\;\;
\vy_{0}^{(2,n-m-1)}\;\;\ldots\;\;\vy_0^{(m+1,n-m-1)}, -\vy_0^{(0,
n-m-1)}),
$$
By assumption, all column vectors of $H$ together with column
vectors of $K^{-}$ are linearly independent. This fact implies that
the left $(m+2)m$ columns of $\vA^*$ together with $(k-1)m$ right
columns of $\vB$ are linearly independent, hence $\rk (\vM) \ge
m(m+k+1)$.  Finally, $\dim \SSS\le (k+m)(m+1) - m(m+k+1)= k$.

We finish the proof by observing that the span of $\vy_0, \vy_1,
\ldots, \vy_{k-1}$ belongs to $\RRR$ and its dimension is $k$
therefore these to subspaces coincide. The vector $\vy_k$ which is
linearly independent with $\vy_0, \vy_1, \ldots, \vy_{k-1}$ can not
belong to $\RRR$ and therefore one of its components $\vy_k^{n-m}$
is linearly independent together with $\SSS$.
\endproof

In this paper, we will apply Proposition~\ref{prop1} for the case
$k=1$ and the pairs of consecutive minimal points. We formulate it
in the corollary.

\begin{corollary}\label{corl_prop1}
Let $m$ be an integer such that $1\le m\le n/2 -1$. Assume that $\lambda
> (n-m)^{-1}$. Then for any two consecutive minimal points
$\vx_{i-1}, \vx_i$ of $\xi$ the vectors $\vx_i^{(0,n-m)}, \ldots,
\vx_i^{(m,n-m)}$ are linearly independent. On top of that, at least
one of the vectors $\vx_{i-1}^{(0,n-m)}, \ldots,
\vx_{i-1}^{(m,n-m)}$ is linearly independent with $\vx_i^{(0,n-m)},
\ldots, \vx_i^{(m,n-m)}$.
\end{corollary}

The next lemma uses the ideas of Laurent~\cite{laurent_2003} and
allows to transfer the linear relation $\va\cdot \vx_i=0$ from a
minimal point $\vx_i$ to the next (or previous) one $\vx_{i\pm 1}$.

\begin{proposition}\label{prop3}
Let $\va\in \langle \vx_i\rangle^\bot \cap \ZZ^{n+1}$. There exists
an absolute constant $c>0$ which satisfies the following conditions.
\begin{itemize}
\item If $||\va||\le \frac{cX_i}{X_{i+1}L_i}$ then $\va\in \langle x_{i+1}\rangle^\bot$.
\item If $||\va||\le cL_{i-1}^{-1}$ then $\va\in \langle x_{i-1}\rangle^\bot$.
\end{itemize}
\end{proposition}

\proof For $\va = (a_0, a_1,\ldots, a_n)$, consider the polynomial
$P_\va(\xi):= \sum_{j=0}^n a_j\xi^j$. In view of $\va\cdot \vx_i=0$,
one computes:
\begin{equation}\label{eq1}
|x_{i,0}P_\va(\xi)| = |\va\cdot \vx_i + a_1(x_{i,0}\xi - x_{i,1}) +
a_2(x_{i,0}\xi^2 - x_{i,2}) + \cdots + a_n(x_{i, 0}\xi^n -
x_{i,n})|\ll ||\va||L_i.
\end{equation}
Recall that we have $x_{i,0}\asymp X_i$. Therefore $|P_\va(\xi)|\ll
\frac{||\va||L_i}{X_i}$.

To prove the first statement, we do similar computations for
$x_{i+1,0}P_\va(\xi)$:
\begin{equation}\label{eq22}
|\va\cdot \vx_{i+1}| \le |x_{i+1,0} P_\va(\xi)| + |a_1(x_{i+1,0}\xi
- x_{i+1,1}) + \cdots + a_n(x_{i+1,0}\xi^n - x_{i+1,n})|
\end{equation}
$$ \ll ||\va||L_{i+1} + \frac{||\va||X_{i+1}L_i}{X_i} \ll
\frac{||\va||X_{i+1}L_i}{X_i}.
$$
Hence, there exists an absolute constant $c>0$ such that if
$||\va||\le \frac{c X_i}{X_{i+1}L_i}$ then $|\va\cdot \vx_{i+1}| <
1$ which in turn implies $\va\in \langle x_{i+1}\rangle^\bot$.

For the second statement, we do the same computation but with
$\va\cdot \vx_{i-1}$:
$$
|\va\cdot \vx_{i-1}| \le |x_{i-1,0} P_\va(\xi)| + |a_1(x_{i-1,0}\xi
- x_{i-1,1}) + \cdots + a_n(x_{i-1,0}\xi^n - x_{i-1,n})|
$$$$
\ll ||\va||L_{i-1} + \frac{||\va||X_{i-1}L_i}{X_i} \ll
||\va||L_{i-1}.
$$
Therefore, for small enough absolute constant $c>0$, if $||\va||\le
c L_{i-1}^{-1}$ then $\va\cdot\vx_{i-1}=0$ and $\va\in \langle
x_{i-1}\rangle^\bot$.
\endproof

We end this section with the result for the case $n=3$. It is
essentially an adaptation of Proposition~5.2 from~\cite{roy_2008}.

\begin{lemma}\label{lem1}
Assume that $n=3$, $\xi$ is transcendental and $\lambda>\sqrt{2}-1$.
Then for all large enough $i\in\NN$, at least one of the numbers
$|\vx_{i-1}^{(0,2)}\wedge \vx_i^{(0,2)}\wedge \vx_i^{(1,2)}|$ or
$|\vx_{i-1}^{(1,2)}\wedge \vx_i^{(0,2)}\wedge \vx_i^{(1,2)}|$ is
non-zero. That immediately implies $X_i\gg
L_i^{\frac{1}{\lambda-1}}\gg X_{i+1}^{\frac{\lambda}{1-\lambda}}$.
\end{lemma}

\proof Proposition~5.2 from~\cite{roy_2008} implies that for
$\lambda>\sqrt{2}-1$ the subspaces $\langle \vx_i^{(0,2)},
\vx_i^{(1,2)}\rangle$ and $\langle \vx_{i-1}^{(0,2)},
\vx_{i-1}^{(1,2)}\rangle$ do not coincide for all large $i$. By
Lemma~2.3 from~\cite{roy_2008} we have that $\vx_i^{(0,2)}$ and
$\vx_i^{(1,2)}$ are linearly independent for all large $i$.
Therefore one of the numbers $|\vx_{i-1}^{(0,2)}\wedge
\vx_i^{(0,2)}\wedge \vx_i^{(1,2)}|$ and $|\vx_{i-1}^{(1,2)}\wedge
\vx_i^{(0,2)}\wedge \vx_i^{(1,2)}|$ is nonzero. By
Proposition~\ref{prop2}, both numbers are estimated as
$$1\le \max\{|\vx_{i-1}^{(0,2)}\wedge
\vx_i^{(0,2)}\wedge \vx_i^{(1,2)}|, |\vx_{i-1}^{(1,2)}\wedge
\vx_i^{(0,2)}\wedge \vx_i^{(1,2)}|\} \ll X_iL_iL_{i-1}\ll
X_i^{1-\lambda} L_i.
$$

Thus we have $X_i\gg L_i^{\frac{1}{\lambda-1}}$. The second
inequality follows from $L_i\ll X_{i+1}^{-\lambda}$.
\endproof

\section{Unconditional improvement of an upper bound for
$\lambda$}

We now prove Theorem~\ref{th3}. Let $m\in\NN$, $n=2m+2$ and
$\lambda>(m+2)^{-1}$. Then Lemma~\ref{lem12} implies that the
vectors $(\vx_i^{(j,n-m-1)})_{0\le j\le m+1}$ are linearly
independent for all large values of $i$. That gives us the estimate
\begin{equation}\label{eq36}
1 \le \left|\left| \bigwedge_{j=0}^{m+1}
\vx_i^{(j,n-m-1)}\right|\right| \stackrel{\eqref{prop2_eq1}}\ll
X_iL_i^{m+1}\quad \stackrel{\eqref{eq24}}\Longrightarrow \quad
X_i\gg X_{i+1}^{(m+1)\lambda}.
\end{equation}

Moreover, it follows from Corollary~\ref{corl_prop1} that at least
one of the vectors $\vx_{i-1}^{(l,n-m)}$ for $0\le l\le m$ is
linearly independent with the system of vectors
$V_i:=(\vx_i^{(j,n-m)})_{0\le j\le m}$.

Suppose that there exist arbitrarily large $i$ such that there are
two vectors $\vx_{i-1}^{(l_1, n-m)}$, $\vx_{i-1}^{(l_2, n-m)}$ which
are linearly independent together with $V_i$. In this case we have a
stronger condition
$$
1 \le \left|\left| \bigwedge_{j=0}^m \vx_i^{(j,n-m)}\wedge
\vx_{i-1}^{(l_1,n-m)}\wedge \vx_{i-1}^{(l_2,n-m)}\right|\right|
\stackrel{\eqref{prop2_eq1}}\ll X_iL_i^mL_{i-1}^2
$$
and therefore, by~\eqref{eq24},
$$
X_i \gg X_i^{2\lambda} X_{i+1}^{m\lambda}.
$$
Since $\lambda>(m+2)^{-1}$, for large enough $i$ we get
$X_i>X_{i+1}$ which is impossible.  We conclude that the assumption
is not satisfied and therefore the dimension of any space $\langle
V_i, V_{i+1}\rangle$ equals $m+2$.

Note that there exist arbitrarily large values of $i$ such that two
spaces $\langle V_{i-1}, V_i\rangle$ and $\langle V_i,
V_{i+1}\rangle$ do not coincide. Indeed, otherwise for large enough
values of $i$, all $\vx_i$ lie in a proper subspace of $\ZZ^{n+1}$
and hence $\xi$ is algebraic. Therefore there exist arbitrarily
large $i$ and the values $0\le k,l\le m$ such that
$\vx_{i-1}^{(k,n-m)}, \vx_i^{(l,n-m)}$ together with $V_{i+1}$ are
linearly independent. Proposition~\ref{prop2} then implies
$$
1 \le \left|\left| \bigwedge_{j=0}^m \vx_{i+1}^{(j,n-m)}\wedge
\vx_{i-1}^{(k,n-m)}\wedge \vx_{i}^{(l,n-m)}\right|\right|
\stackrel{\eqref{prop2_eq1}}\ll X_{i+1}L_{i+1}^mL_{i-1}L_i
\stackrel{\eqref{eq24}}\ll X_{i+1}^{1-\lambda} X_{i+2}^{-m\lambda}
X_i^{-\lambda}.
$$
By~\eqref{eq36}, we further bound it from above by
$$
X_{i+1}^{ 1-(m+1)\lambda - ( m+1)\lambda^2}.
$$
Therefore if $\lambda$ is bigger than the positive root $\lambda_0$
of $(m+1)x^2+(m+1)x - 1=0$, then this expression becomes arbitrarily
close to zero which is a contradiction. Observe that
$(m+2)^{-1}<\lambda_0<(m+1)^{-1}$ which finally gives that
$\lambda\le \lambda_0$. The proof of Theorem~\ref{th3} is finished.

\section{Relations between $\lambda$ and $\omega_k$}

Given a positive parameter $Y$ and a positive integer $k$, consider
a centrally-symmetric convex figure $F(Y)\subset \RR^{k+1}$ defined
by the inequalities
$$
F(Y):=\{\vy\in \RR^{k+1}\;:\; ||\vy^{(1,k-1)}||\le Y;\ |y_0 + y_1\xi + \ldots + y_k\xi^k|\le Y^{-k}\}.
$$
We first establish an estimate on the last successive minimum of $F(Y)$,
based on $\omega_k$.

\begin{lemma}\label{lem13}
For all $Y\ge 1$, the $(k+1)$'th successive minimum $\tau_{k+1}$ of
$F(Y)$ is bounded from above by
\begin{equation}\label{lem13_eq1}
\tau_{k+1} \ll Y^{\frac{k(\omega_k-k)}{1+\omega_k}}
\end{equation}
\end{lemma}

\proof We first note that the first successive minimum $\tau_1$ of
$F(Y)$ is bigger than $cY^{\frac{k-\omega_k}{1+\omega_k}}$ for some
small absolute constant. Indeed, if the inequalities
$$
||\vy||\le Y^{\frac{k-\omega_k}{1+\omega_k}}\cdot Y =
Y^{\frac{k+1}{1+\omega_k}};\quad |y_0 + y_1\xi + \cdots +
y_k\xi^k|\le Y^{\frac{k-\omega_k}{1+\omega_k}} Y^{-k} =
Y^{-\frac{(k+1)\omega_k}{1+\omega_k}}
$$
are satisfied then we get that the inequalities~\eqref{def_omk} have
an integer solution $\vy$ for $X = Y^{\frac{k+1}{1+\omega_k}}$. If
$X$ lies outside a bounded region from~\eqref{def_omk} then we
immediately get a contradiction. For smaller values of $X$ we may
have a finite number of solutions $\vy$ of the inequalities above.
All of them can be ruled out by taking the constant $c$ in the
expression $\tau_1\ge cY^{\frac{k-\omega_k}{1+\omega_k}}$ small
enough.

Secondly, the volume of $F(Y)$ is $2^{k+1}$ and hence by the
Minkowski's theorem on successive minima, we have
$\tau_1\tau_2\cdots\tau_{k+1}\asymp 1$. That gives an upper bound
$$
\tau_{k+1} \ll \tau_1^{-k} \ll Y^{\frac{k(\omega_k-k)}{1+\omega_k}}.
$$
\endproof

For a given vector $\va= (a_0, a_1,\ldots, a_k)\in\RR^{k+1}$ we
denote by $P_\va(x)$ a polynomial $a_0 + a_1x+\cdots + a_kx^k$. And
we also define an inverse operation: for a given a polynomial
$P\in\RR[x]$ of degree $k$ we denote by $\va(P)$ a vector in
$\RR^{k+1}$ composed from the polynomials coefficients.

\begin{lemma}\label{lem14}
let $P\in \RR[x]$ be a polynomial of degree $m$ and
$\vx\in\RR^{m+k+1}$. Suppose there exist $k+1$ linearly independent
polynomials $Q_i\in \RR[x]$, $i\in\{0,\ldots, k\}$ of degree at most
$k$ such that for all $0\le i\le k$, $\va(Q_iP)\in \langle
\vx\rangle^\bot$. Then one has $\va(P)\in \langle \vx^{(0,m)},
\vx^{(1,m)}, \ldots, \vx^{(k,m)}\rangle^\bot$.
\end{lemma}

\proof The space of polynomials of degree at most $k$ has dimension
$k+1$, therefore the polynomials $Q_i$ form its basis. That in turn
implies that for each $j\in\{0,\ldots, k\}$ one can write $x^j$ as a
linear combination
$$
x^j = \sum_{i=0}^k c_{j,i} Q_i(x).
$$
We then have $\va(x^j\cdot P(x))\in \langle \vx\rangle^\bot$ which
is equivalent to $\va(P)\in \langle \vx^{(j,m)}\rangle^\bot$. To
finish the proof, we observe that
$$
\langle \vx^{(0,m)}, \vx^{(1,m)}, \ldots, \vx^{(k,m)}\rangle^\bot =
\bigcap_{i=0}^k \langle \vx^{(i,m)}\rangle^\bot.
$$
\endproof

Fix a positive integer $k$. The key idea of the proof below is to
consider the shortest non-zero vector $\va$ from $\langle
\vx_i^{(0,n-m)}, \ldots, \vx_i^{(m,n-m)}\rangle^\bot \cap
\ZZ^{n-m+1}$ where $m\ge k$ and construct $k+1$ linearly independent
polynomials $Q_j$, $0\le j\le k$ such that $Q_jP_\va$ belongs to
$\langle \vx_{i+1}\rangle^\bot$. That will imply that $\va$ belongs
to $\langle \vx_{i+1}^{(0,n-m)}, \ldots,
\vx_{i+1}^{(m,n-m)}\rangle^\bot$ and by repeating this process again
and again we infer that vectors $\vx_i$ for all large $i$ belong to
a proper subspace of $\RR^{n+1}$ which contradicts to that $\xi$ is
transcendental.

Denote $\Lambda_i:= \langle \vx_{i}^{(0,n-m)}, \ldots,
\vx_{i}^{(m,n-m)}\rangle^\bot\cap \ZZ^{n-m+1}$. For the rest of this
section, let $\va = \va_i$ be the shortest non-zero vector in
$\Lambda_i$. The latter set is a sublattice in $\RR^{n-m+1}$. By
Lemma~\ref{lem12}, if $m\le \frac{n}{2}$ and $\lambda>(n-m+1)^{-1}$
then the vectors $\vx_i^{(j,n-m)}$ are linearly independent and the
dimension of $\Lambda_i$ is $n-2m$. Therefore $\Lambda_i$ is indeed
a non-trivial sublattice and not just a single zero point if $n\ge
2m+1$ or $m\le \lfloor\frac{n-1}{2}\rfloor$.

The covolume of $\Lambda_i$ in the subspace $\langle
\vx_{i}^{(0,n-m)}, \ldots, \vx_{i}^{(m,n-m)}\rangle^\bot$ equals the
length of the primitive multi-vector which is a scalar multiple of
$$ \bigwedge_{j=0}^m \vx_i^{(j,n-m)}. $$
Therefore by the Minkowski's theorem and Proposition~\ref{prop2},
the shortest non-zero vector of $\Lambda_i$ satisfies
\begin{equation}\label{eq18}
||\va||\ll \left|\left|\bigwedge_{j=0}^m
\vx_i^{(j,n-m)}\right|\right|^{\frac{1}{n-2m}} \ll (X_i
L_i^m)^{\frac{1}{n-2m}}.
\end{equation}
The linear independence of vectors~$\vx_i^{(j,n-m)}$ then gives us
that $X_iL_i^m\gg 1$. Hence by~\eqref{eq24}, we get that for all
large $i$,
\begin{equation}\label{eq26}
X_i\gg X_{i+1}^{m\lambda}.
\end{equation}

Corollary~\ref{corl_prop1} may give us a better estimate on
$||\va||$ in some cases. If $\lambda> (n-m)^{-1}$ then one of the
vectors $\vx_{i-1}^{(l,n-m)}$, $0\le l\le m$ is linearly independent
together with $(\vx_i^{(j,n-m)})_{0\le j\le m}$. Therefore the
shortest non-zero vector of $\Lambda_i$ is not longer than the
shortest non-zero vector of $\Lambda_i\cap \langle
\vx_{i-1}^{(l,n-m)}\rangle^\bot$. The latter is a sublattice of
dimension $n-2m-1$ which contains non-zero vectors for $m\le
\lfloor\frac{n-2}{2}\rfloor$. The application of Minkowski theorem
and Proposition~\ref{prop2} leads to
\begin{equation}\label{eq19}
||\va||\ll \left|\left|\bigwedge_{j=0}^m \vx_{i}^{(j,n-m)}\wedge
\vx_{i-1}^{(l,n-m)}\right|\right|^{\frac{1}{n-2m-1}} \ll (X_i
L_i^mL_{i-1})^{\frac{1}{n-2m-1}}.
\end{equation}

Conditions $\lambda>(n-m)^{-1}$ together with $n\ge 2m+2$ imply that
the vectors $\vx_i^{(j,n-m-1)}$, $j\in \{0,\ldots, m+1\}$ are
linearly independent which, analogously to~\eqref{eq26} implies that
\begin{equation}\label{eq25}
X_i\gg X_{i+1}^{(m+1)\lambda}.
\end{equation}

We denote the length of $\va$ by $A$ and keep in mind that both
bounds~\eqref{eq18} and~\eqref{eq19}, with slightly different
conditions on $m$ and $\lambda$, are satisfied.

Note that
$$
|x_{i,0}P_\va(\xi)| = |a_1(x_{i,0}\xi - x_{i,1}) + \cdots +
a_m(x_{i,0}\xi^m - x_{i,m})| \ll AL_i.
$$
Therefore $|P_\va(\xi)|\ll \frac{AL_i}{X_i}$.

By Lemma~\ref{lem13}, for any large enough $Y$ there exist $k+1$ linearly
independent polynomials $Q_j\in \ZZ[x]$, $0\le j\le k$ such that
\begin{equation}\label{eq20}
\deg Q_j\le k,\ ||\va(Q_j)||\ll
Y^{\frac{k(\omega_k-k)}{1+\omega_k}+1},\ |Q_j(\xi)|\ll
Y^{\frac{k(\omega_k-k)}{1+\omega_k} - k}.
\end{equation}
Choose $Y$ such that for all $j\in\{0,\ldots, k\}$ one has $|x_{i+1,0}
Q_j(\xi)P_\va(\xi)| < \frac12$. We have that
$$
|x_{i+1,0} Q_j(\xi)P_\va(\xi)| \ll \frac{X_{i+1}AL_i}{X_i}
Y^{-\frac{k(k+1)}{1+\omega_k}}.
$$
Hence,
\begin{equation}\label{eq21}
Y \ge C
\left(\frac{X_{i+1}AL_i}{X_i}\right)^{\frac{1+\omega_k}{k(k+1)}}
\asymp
\left(\frac{X_{i+1}AL_i}{X_i}\right)^{\frac{1+\omega_k}{k(k+1)}}
\end{equation}
does the job, assuming that $C$ is a large enough absolute constant.

Let $\vb = \vb(j) := \va(Q_jP_\va)\in \ZZ^{n-m+k+1}$. Similarly
to~\eqref{eq22} we compute
\begin{equation}\label{eq23}
|\vb\cdot \vx_{i+1}^{(0,n-m+k)}| \le |x_{i+1,0}Q_jP_\va(\xi)| +
|b_1(x_{i+1,0}\xi-x_{i+1,1}) + \cdots + b_{n-m+k}(x_{i+1,0}\xi^{n-m+k} -
x_{i+1,n-m+k})|
\end{equation}
The second term on the right hand side is $\ll ||\vb||L_{i+1}$
which, by~\eqref{eq20}, is bounded from above by
\begin{equation}\label{eq27}
Y^{\frac{(k+1)(\omega_k-k+1)}{1+\omega_k}} AL_{i+1}\ll
\left(\frac{X_{i+1}AL_i}{X_i}\right)^{\frac{\omega_k-k+1}{k}}\cdot
AL_{i+1} =
\left(\frac{X_{i+1}L_i}{X_i}\right)^{\frac{\omega_k-k+1}{k}}A^{\frac{\omega_k+1}{k}}
L_{i+1}.
\end{equation}
Now, if the last expression is smaller than some small enough
absolute constant, so that the second term in~\eqref{eq23} is
smaller than $\frac12$ then we derive that $\vb(j)\cdot
\vx_{i+1}^{(0,n-m+k)}=0$ for all $j\in\{0,\ldots, k\}$. By analogy,
under the same conditions we derive that $\vb(j)\cdot
\vx_{i+1}^{(l,n-m+k)}=0$ for all $j\in \{0,\ldots, k\}$ and
$l\in\{0,\ldots, m-k\}$. Then Lemma~\ref{lem14} implies that $\va\in
\Lambda_{i+1}$.

We apply the upper bound~\eqref{eq18} for $A$. Recall that it is
satisfied in the case $n\ge 2m+1$ and $\lambda>(n-m+1)^{-1}$. Then
we continue estimating~\eqref{eq27}:
$$
\left(\frac{X_{i+1}L_i}{X_i}\right)^{\frac{\omega_k-k+1}{k}}A^{\frac{\omega_k+1}{k}}
L_{i+1} \ll
X_{i+1}^{\frac{\omega_k-k+1}{k}}X_i^{\frac{\omega_k+1}{k(n-2m)}-
\frac{\omega_k-k+1}{k}}L_i^{\frac{\omega_k-k+1}{k}+\frac{m(\omega_k+1)}{k(n-2m)}}L_{i+1}
$$
\begin{equation}\label{eq28}
\stackrel{\eqref{eq24}}\ll X_{i+1}^{\frac{\omega_k-k+1}{k} -
\frac{(\omega_k+1)(n-m)}{k(n-2m)}\lambda} X_i^{1 -
\frac{(\omega_k+1)(n-2m-1)}{k(n-2m)}} \frac{L_{i+1}}{L_i}.
\end{equation}

We consider two cases.

{\bf Case 1.} The degree of $X_i$ is nonnegative. That is equivalent
to
$$
1\ge \frac{(\omega_k+1)(n-2m-1)}{k(n-2m)}\quad
\Longleftrightarrow\quad n\le 2m+1 + \frac{k}{\omega_k+1-k}.
$$
In this case we can use a straightforward bound $X_i< X_{i+1}$. We
also use $L_{i+1}< L_i$ to finally get an upper bound
for~\eqref{eq28}:
$$
X_{i+1}^{\frac{\omega_k-k+1}{k} -
\frac{(\omega_k+1)(n-m)}{k(n-2m)}\lambda + 1 -
\frac{(\omega_k+1)(n-2m-1)}{k(n-2m)}} =
X_{i+1}^{\frac{\omega_k+1}{k(n-2m)} -
\frac{(\omega_k+1)(n-m)}{k(n-2m)}\lambda}.
$$
Note that for $\lambda>(n-m)^{-1}$ the last expression becomes
arbitrary close to zero as $X_{i+1}\to\infty$. Therefore for large
enough $i$ the shortest vector of $\Lambda_i$ also belongs to
$\Lambda_{i+1}$. This implies that the sequence of lengths
$||\va_i||$ of the shortest non-zero vectors of $\Lambda_i$ is
monotonically non-increasing. But this sequence can not decrease
infinitely often, therefore there exists $i_0\in\NN$ such that for
all $i>i_0$ the sequence $||\va_i||$ is constant, hence all the
vectors $\vx_i$ belong to a proper subspace and thus $\xi$ is
algebraic.

We derive that for $n$ in the range
\begin{equation}\label{eq31}
2m+1\le n\le
2m+1+\frac{k}{\omega_k+1-k} = 2m+1+\delta_k,
\end{equation}
the parameter $\lambda$ must satisfy $\lambda\le (n-m)^{-1}$.

The expression $(n-m)^{-1}$ grows together with $m$, therefore to make an
upper bound on $\lambda$ as small as possible, we need to take the smallest
$m$ such that the condition~\eqref{eq31} is satisfied. In view of the
condition $m\ge k$, if $2k+1\le n\le 2k+1+\delta_k$ then the smallest
possible $m$ is $m=k$ and $\lambda \le (n-k)^{-1}$. Since $\omega_k$ can be
taken arbitrary close to $\omega _k(\xi)$, the last assertion implies the
first statement of Theorem~\ref{th1}.

If $n>2k+1+\delta_k$ then $m\ge \frac{n-\delta_k-1}{2}>k$ i.e. the
smallest possible $m$ is $\lceil \frac{n-\delta_k-1}{2}\rceil$.
Because of the assumption~\eqref{th1_eq1} that $\delta_k\ge 1$, we
always get $2m+1\le n$ and the condition~\eqref{eq31} is satisfied.
Now the inequality $\lambda\le (n-m)^{-1}$ implies the first part of
the second assertion of Theorem~\ref{th1}.

{\bf Case 2.} The degree of $X_i$ is negative which is equivalent to
$$
n> 2m+1 + \frac{k}{\omega_k+1-k}.
$$
In this case we use the bound~\eqref{eq26} to further
estimate~\eqref{eq28}:
$$
X_{i+1}^{\frac{\omega_k-k+1}{k} -
\frac{(\omega_k+1)(n-m)}{k(n-2m)}\lambda + m\lambda  -
\frac{(\omega_k+1)(n-2m-1)}{k(n-2m)}m\lambda} =
X_{i+1}^{\frac{\omega_k-k+1}{k} - \frac{\omega_k-k+1}{k}m\lambda -
\frac{\omega_k+1}{k}\lambda}.
$$
This expression becomes arbitrarily small for $i$ large enough if
$$
\lambda > \frac{\omega_k-k+1}{(\omega_k-k+1)m + \omega_k+1} =
\frac{1}{m+1+\delta_k}.
$$
As in the case~1, for such $\lambda$ we get that the sequence
$||\va_i||$ for the shortest vectors $\va_i$ of lattices $\Lambda_i$
is eventually constant and hence $\xi$ is algebraic.

The last bound on $\lambda$ decays as $m$ grows, therefore the
smallest possible upper bound on $\lambda$ will be when $m$ is the
largest possible such that $2m+1 + \delta_k<n$  or $ m =
\lfloor\frac{n-\delta_k-1}{2}\rfloor$. Taking into account that all
the arguments only work for $\lambda>(n-m+1)^{-1}$, we finally
derive the upper bound:
$$
\lambda\le \max\left\{\frac{1}{n-m+1}, \frac{1}{\left\lfloor
\frac{n-\delta_k-1}{2}\right\rfloor +1+\delta_k}\right\} =
\frac{1}{\left\lfloor \frac{n-\delta_k-1}{2}\right\rfloor
+1+\delta_k}
$$
which is the second term in the second assertion of
Theorem~\ref{th1}.

Finally, we need to check the condition $m\ge k$. It is only
satisfied when $n>2k+1+\delta_k$. That finishes the proof of
Theorem~\ref{th1}.

\section{Light improvement of Theorem~\ref{th1}}\label{sec6}

For some values of $n$ and $\delta_k$, the upper bound on
$\widehat{\lambda}_n(\xi)$ in Theorem~\ref{th1} can be slightly
improved if we use estimates~\eqref{eq19} and~\eqref{eq25} instead
of~\eqref{eq18} and~\eqref{eq26} respectively. However the
corresponding expressions become much more technical. Also, to use
those estimates, we need to make stronger assumptions on $m$ and
$\lambda$: $n\ge 2m+2$ and $\lambda>(n-m)^{-1}$.

\begin{theorem}\label{th2}
Let $\xi\in\RR$ be a transcendental number, $k\in\NN$ and $\delta_k$
be defined by~\eqref{th1_eq1}. Then
$$
\widehat\lambda_n(\xi)\le \frac{1}{n-m}
$$
for the minimal value of $m$ such that $m\ge k$ and
\begin{equation}\label{th2_eq1}
2m+2\le n< 2m+1+\left(1-\frac{1}{n-m}\right) (1+\delta_k).
\end{equation}

For a given $m\in\NN$, consider the positive root $x = x(m)$ of the
quadratic equation
\begin{equation}\label{th2_eq2}
\frac{1}{\delta_k} - \frac{(n-m-1)x}{\delta_k(n-2m-1)} -
\frac{(n-m-1)x}{n-2m-1} + \left(\frac{1-x}{n-2m-1} -
\frac{n-2m-2+x}{\delta_k(n-2m-1)}\right)(m+1)x.
\end{equation}
Then for any $m$ such that $2m+2\le n$, $2m+1+(1-x(m))(1+\delta_k)\le n$ and
$k\le m$ one has $\widehat{\lambda}_n(\xi)\le x(m)$.
\end{theorem}

\proof We apply an upper bound~\eqref{eq19} on $A$ in
estimate~\eqref{eq27}:
$$
\left(\frac{X_{i+1}L_i}{X_i}\right)^{\frac{\omega_k-k+1}{k}}A^{\frac{\omega_k+1}{k}}
L_{i+1} \ll X_{i+1}^{\frac{\omega_k-k+1}{k}}
X_i^{\frac{\omega_k+1}{k(n-2m-1)} - \frac{\omega_k-k+1}{k}}
L_i^{\frac{\omega_k-k+1}{k} +
\frac{m(\omega_k+1)}{k(n-2m-1)}}L_{i-1}^{\frac{\omega_k+1}{k(n-2m-1)}}L_{i+1}
$$
\begin{equation}\label{eq32}
\stackrel{\eqref{eq24}}\ll X_{i+1}^{\frac{\omega_k-k+1}{k} -
\frac{(\omega_k+1)(n-m-1)}{k(n-2m-1)}\lambda}
X_i^{1-\frac{(\omega_k+1)(n-2m-2)}{k(n-2m-1)} -
\frac{(\omega_k+1)}{k(n-2m-1)}\lambda} \frac{L_{i+1}}{L_i}.
\end{equation}

As in the proof of Theorem~\ref{th1}, we consider two cases.

{\bf Case 1*.} The degree of $X_i$ is nonnegative. That is
equivalent to
$$
1\ge \frac{(\omega_k+1)(n-2m-2+\lambda)}{k(n-2m-1)} \quad
\Longleftrightarrow\quad n\le 2m+1+(1-\lambda)(1+\delta_k).
$$
Then the inequalities $X_i<X_{i+1}$ and $L_{i+1}<L_i$ give an upper
bound for~\eqref{eq32}:
$$
X_{i+1}^{\frac{\omega_k+1}{k} -
\frac{(\omega_k+1)(n-m-1)}{k(n-2m-1)}\lambda -
\frac{(\omega_k+1)(n-2m-2)}{k(n-2m-1)} -
\frac{\omega_k+1}{k(n-2m-1)}\lambda} = X_{i+1}^{\frac{(\omega_k+1)(1
- (n-m)\lambda)}{k(n-2m-1)}}.
$$
If $\lambda>(n-m)^{-1}$ we immediately get that the last expression becomes
arbitrary close to zero as $X_{i+1}\to\infty$. By analogy with Case~1, that
implies that the sequence of $\va_i$ is constant for $i>i_0$ and $\xi$ is
algebraic.

If $n$ is such that
\begin{equation}\label{eq33}
2m+2\le n< 2m+1 +\left(1 - \frac{1}{n-m}\right) (1+\delta_k)
\end{equation}
then there exists $\lambda = (n-m)^{-1} + \epsilon$ for small enough
$\epsilon$ such that $n\le 2m+1 + (1-\lambda)(1+\delta_k)$ and we are in Case
1*. However, if the upper bound~\eqref{eq33} on $n$ is not satisfied then the
degree of $X_i$ in~\eqref{eq32}  is never positive for $\lambda>(n-m)^{-1}$
and hence Case~1* does not take place. This give us the first assertion of
Theorem~\ref{th2}.

{\bf Case 2*.} The degree of $X_i$ is negative which is equivalent
to $n>2m+1+(1-\lambda)(1+\delta_k)$. Then we use~\eqref{eq25} to
further estimate~\eqref{eq32}:
$$
\ll X_{i+1}^{\frac{\omega_k-k+1}{k} -
\frac{(\omega_k+1)(n-m-1)}{k(n-2m-1)}\lambda +
\left(1-\frac{(\omega_k+1)(n-2m-2)}{k(n-2m-1)} -
\frac{\omega_k+1}{k(n-2m-1)}\lambda\right)(m+1)\lambda}
$$
We slightly simplify the degree of $X_{i+1}$:
$$
\frac{\omega_k-k+1}{k} -
\frac{(\omega_k+1)(n-m-1)}{k(n-2m-1)}\lambda +
\left(1-\frac{(\omega_k+1)(n-2m-2)}{k(n-2m-1)} -
\frac{(\omega_k+1)\lambda}{k(n-2m-1)}\right)(m+1)\lambda
$$
\begin{equation}\label{eq34}
= \frac{1}{\delta_k} - \frac{(n-m-1)\lambda}{\delta_k(n-2m-1)} -
\frac{(n-m-1)\lambda}{n-2m-1} + \left(\frac{1-\lambda}{n-2m-1} -
\frac{n-2m-2+\lambda}{\delta_k(n-2m-1)}\right)(m+1)\lambda.
\end{equation}
Note that this expression is exactly~\eqref{th2_eq2} with $\lambda$
instead of $x(m)$. For $\lambda = (n-m)^{-1}$ one has $(m+1)\lambda
<1$ and hence the expression~\eqref{eq34}, as in Case~1*, is bounded
from below by
$$
\frac{(\omega_k+1)(1-(n-m)\lambda)}{k(n-2m-1)} = 0.
$$
Therefore we have $x(m) > (n-m)^{-1}$ and any $\lambda > x(m)$ also
satisfies the condition $\lambda> (n-m)^{-1}$. While it is not
needed for the proof itself, it is useful to observe that by the
same method we can check that $x(m)$ lies between $(n-m)^{-1}$ and
$(m+1)^{-1}$.

For $\lambda > x(m)$ the expression~\eqref{eq34} becomes negative
and the upper bound in~\eqref{eq32} becomes arbitrary close to zero
as $X_{i+1}\to\infty$ and, as before, this implies that $\xi$ is
algebraic. Therefore we must have $\lambda\le x(m)$. \endproof

To demonstrate that Theorem~\ref{th2} indeed gives slightly better
estimates in some cases, we consider a couple of examples. Observe
that, by Dirichlet theorem, $\omega_k$ can take values not smaller
than $k$ and therefore $\delta_k$ takes values between 0 and $k$.

\begin{example}
$k=1$ and $\delta_k$=1. That is, $\xi$ is not very well approximable.
\end{example}

We emphasize that the conditions of Example~1 are satisfied for
almost all real numbers $\xi$ in terms of the Lebesgue measure.

In this case, Theorem~\ref{th1} gives us that for any $n\ge 2$,
$$
\widehat\lambda_{2n}(\xi) \le \frac{1}{n+1},\quad \widehat\lambda_{2n-1}(\xi)\le \frac{1}{n}.
$$
The achieved upper bound on $\widehat\lambda_{2n}(\xi)$ is currently
best known (of course, it depends on the condition $\delta_1=1$).
However the bound on $\widehat\lambda_{2n-1}(\xi)$ coincides with
that of Laurent~\cite{laurent_2003}.

Now we apply Theorem~\ref{th2}. Notice that for odd $n$ the
condition~\eqref{th2_eq1} is not satisfied for all integer values of
$m$, hence the first assertion of this theorem can not be applied.
On the other hand, for $n\ge 2$ this assertion gives
$$
\widehat\lambda_{2n}(\xi)\le \frac{1}{n+1}.
$$
That is exactly the same bound as in Theorem~\ref{th1}. Now we apply the
second assertion of Theorem~\ref{th2}. The equation~\eqref{th2_eq2} in this
case simplifies to
$$
2(m+1)x^2 + (m(n-2m) + 3n - 7m-5)x - (n-2m-1)=0.
$$
If $n = 2m+3$ then the equation simplifies even further to $(m+1)x^2
+ (m+2)x - 1=0$. One can check that its positive root $x(m)$ belongs
to $\big(\frac{1}{m+3}, \frac{1}{m+2}\big)$. Simple calculations
then verify that all the conditions of Theorem~\ref{th2} are
satisfied and
$$
\widehat{\lambda}_{2m+3}\le x(m) = \frac{\sqrt{m^2 +
8m+8}-(m+2)}{2(m+1)}.
$$
This upper bound is better than $\frac{1}{m+2}$ derived from
Theorem~\ref{th1}. For example,
$$
\widehat\lambda_5(\xi)\le \frac{\sqrt{17}-3}{4}\approx
0.2808<\frac13;\quad \widehat\lambda_7(\xi)\le
\frac{\sqrt{7}-2}{3}\approx 0.2153 < \frac14.
$$

\begin{example}
$k=2$ and $\delta_k = 3/2$ which in turn means that $\omega_2(\xi) =
7/3$.
\end{example}

Theorem~\ref{th1} gives that
$$
\widehat\lambda_5(\xi)\le \frac13,\quad
\widehat\lambda_{2n}(\xi)\le\frac{1}{n+1},\quad
\widehat\lambda_{2n+1}(\xi)\le \frac{2}{2n+3}\quad \forall n\ge 3.
$$
All these upper bounds are smaller than the best currently known
bounds on $\widehat\lambda_n(\xi)$.

The first assertion of Theorem~\ref{th2} gives the same upper bounds
on $\widehat\lambda_{2n}(\xi)$, $n\ge 3$ as Theorem~\ref{th1}. Also,
the conditions of the first assertion are never satisfied for $n=5$
and $n=7$. However, for $m\ge 3$ it gives
\begin{equation}\label{eq35}
\widehat\lambda_{2m+3}(\xi)\le \frac{1}{m+3},
\end{equation}
which is better than in Theorem~\ref{th1}.

The equation~\eqref{th2_eq2} simplifies to
$$
5(m+1)x^2 + (2m(n-2m) + 7n - 16m-12)x - 2(n-2m-1)=0
$$
For $n=7, m=2$ we get the solution $x(m) = \frac15$ and the second
assertion of Theorem~\ref{th2} infers $\widehat\lambda_7(\xi)\le
\frac15$ which complements~\eqref{eq35} for $m=2$. If $n$ and $m$
are related by $n=2m+4$, the equation further simplifies to
$$
5(m+1)x^2 + (6m+16)x - 6=0.
$$
By substituting $x = (m+3)^{-1}$ and $x = (m+4)^{-1}$ into the
quadratic polynomial above and checking that the result is positive
in the first case and is negative in the second one, one concludes
that $x(m)\in \big(\frac{1}{m+4}, \frac{1}{m+3}\big)$ which is
smaller than the value $\frac{1}{m+3}$ provided by
Theorem~\ref{th1}. In particular,
$$
\widehat\lambda_8(\xi)\le \frac{\sqrt{286}-14}{15}\approx
0.1941<\frac15;\quad \widehat\lambda_{10}(\xi)\le
\frac{\sqrt{409}-17}{20}\approx 0.1612<\frac16.
$$
We can see that in Example~2, Theorem~\ref{th2} gives better upper
bounds for $\widehat\lambda_n(\xi)$ for all values of $n\ge 7$.

\section{Better relation between $\lambda$ and
$\omega_1$ for $n=3$}

Here we prove Theorem~\ref{th3}. Denote by $\vc = \vc_i$ the
shortest vector in the subspace $\langle\vx_i\rangle^\bot \cap
\ZZ^4$. One can easily check that for transcendental $\xi$ there
must exist arbitrarily large values of $i$ such that $\vc_i \not\in
\langle\vx_{i+1}\rangle^\bot \cap \ZZ^4$. Fix one such value of $i$.
By Proposition~\ref{prop3}, we have that
\begin{equation}\label{eq38}
||\vc||\gg \frac{X_i}{X_{i+1}L_i} \gg
\frac{X_i}{X_{i+1}^{1-\lambda}}.
\end{equation}

By the Minkowski's theorem, $||\vc||\le X_i^{1/3}$ which is in turn
smaller than $L_{i-1}^{-1}$. therefore Proposition~\ref{prop3}
implies that $\vc_i$ also belongs to $\langle\vx_{i-1}\rangle^\bot
\cap \ZZ^4$. Since $\vc_i \not\in \langle\vx_{i+1}\rangle^\bot \cap
\ZZ^4$, we get that $\vx_{i-1}, \vx_i$ and $\vx_{i+1}$ are linearly
independent or equivalently, $i\in I$.

Let $j$ be the predecessor of $i$ in the set $I$. Since $I$ is an
infinite set, we can choose $i$ large enough so that $j$ exists.
Then we have $\langle \vx_j, \vx_{j+1}\rangle = \langle \vx_{i-1},
\vx_i\rangle$ and any vector from $\langle \vx_{j-1}, \vx_j,
\vx_{j+1}\rangle^\bot$ is orthogonal to $\vx_i$. Therefore by
Proposition~\ref{prop2} and~\eqref{eq24}, we get
\begin{equation}\label{eq39}
||\vc_i|| \le ||\vx_{j-1}\wedge \vx_j\wedge \vx_{j+1}|| \ll
X_{j+1}^{1-\lambda}X_j^{-\lambda}.
\end{equation}
Combining two inequalities~\eqref{eq38} and~\eqref{eq39} together
gives us
\begin{equation}\label{eq41}
X_{j+1}^{1-\lambda}X_{i+1}^{1-\lambda} \gg X_j^\lambda X_i.
\end{equation}

From now on, we assume that $\lambda>\sqrt{2}-1$ which allows us to
apply Lemma~\ref{lem1}. Since $\vx_i\in\langle \vx_j,
\vx_{j+1}\rangle$ we have $\vx_i = u\vx_j+v\vx_{j+1}$ for some real
$u$ and $v$. Also, since $\vx_i$ is neither a scalar multiple of
$\vx_j$ nor of $\vx_{j+1}$, we have that $u$ and $v$ are non-zero.

The following arguments are adapted from the proof of Lemma~4.2
in~\cite{roy_2008}. By Lemma~3.1 from~\cite{roy_2008}, the
coefficients $u$ and $v$ are integer. If $X_i>3|v|X_{j+1}$ we have
$|u|X_j = ||\vx_i - v\vx_{j+1}|| > 2|v|X_{j+1}$ and so $|u|>2|v|$.
But then we find $L_i\ge |u|L_j - |v|L_{j+1} > L_{j+1}$ which is
impossible. This contradiction shows that $|v|\ge X_i/(3X_{j+1})\gg
\frac{X_i}{X_{j+1}}$.

Now for $l\in\{1,2\}$ we have that
$$
||\vx_j^{(0,2)}\wedge \vx_j^{(1,2)}\wedge \vx_i^{(l,2)}|| = v
||\vx_j^{(0,2)}\wedge \vx_j^{(1,2)}\wedge \vx_{j+1}^{(l,2)}||.
$$
By Lemma~\ref{lem1}, we have that at least one of
$\vx_{j+1}^{(0,2)}$ or $\vx_{j+1}^{(1,2)}$ is linearly independent
with $(\vx_j^{(0,2)},\vx_j^{(1,2)})$ and hence the same is true for
at least one of the vectors $\vx_i^{(0,2)}$ or $\vx_i^{(1,2)}$.

Now we will closely follow the proof of Theorem~\ref{th1} with
$k=m=1$. Let $\va = \va_j$ be the shortest non-zero vector in
$\Lambda_j$. Then~\eqref{eq18} gives us $A:=||\va||\ll X_jL_j$.
Observe that $|x_{j,0}P_\va(\xi)| \ll AL_j$ and hence
$|P_\va(\xi)|\ll \frac{AL_j}{X_j}$.

Choose two linearly independent polynomials $Q_0, Q_1\in\ZZ[x]$ such
that~\eqref{eq20} is satisfied for $k=1$. By Lemma~\ref{lem13}, We
can do that for any positive $Y$. Choose $Y$ such that for both
$l\in\{0,1\}$ one has $|x_{i,0} Q_l(\xi)P_\va(\xi)|<\frac12|v|$.
Then for a small enough absolute constant $c$ the value
\begin{equation}\label{eq40}
Y = c\left(\frac{X_iAL_j}{|v| X_j}\right)^{\frac{1+\omega_1}{2}}
\end{equation}
does the job.

Let $\vb = \vb(l):= \va(Q_lP_\va)\in\ZZ^{n+1}$. We then have
$\vb\cdot \vx_j=0$. Indeed, if $Q_l(\xi) = c_1\xi + c_0$ then
$\vb\cdot \vx_j = (0,c_1\va)\cdot \vx_j + (c_0\va,0)\cdot \vx_j$ and
the last expression equals zero since $\va\in\Lambda_j$.
From~\eqref{eq20} we have that $||\vb||\ll
Y^{\frac{2\omega_1}{1+\omega_1}}A$. Next, we compute
$$
|\vb\cdot \vx_i|\le |x_{i,0} Q_l(\xi)P_\va(\xi)| + |b_1(x_{i,0}\xi -
x_{i,1}) + \cdots + b_3(x_{i,0}\xi^3 - x_{i,3})|.
$$
The left hand side is equal to $|\vb\cdot \vx_i| = |\vb \cdot
(u\vx_j+v\vx_{j+1})| = |v| |\vb\cdot \vx_{j+1}|$ and we derive that
it is a multiple of $v$. If the second term on the right hand side
is smaller than $\frac12|v|$ then we get $\vb(0)\cdot
\vx_i=\vb(1)\cdot \vx_i =0$ which by Lemma~\ref{lem14} implies that
$\va_j\in\Lambda_{i}$. Finally, that means that the systems
$\vx_j^{(0,2)}, \vx_j^{(1,2)},\vx_i^{(0,2)}$ and $\vx_j^{(0,2)},
\vx_j^{(1,2)},\vx_i^{(1,2)}$ are linearly dependent which is
impossible, as we showed above.

We conclude that the second term must be at least $\frac12|v|$
which, together with~\eqref{eq40} and $A\ll X_jL_j$ infers
$$
\frac{X_i}{X_{j+1}}\ll|v|\ll Y^{\frac{2\omega_1}{1+\omega_1}} AL_i
\ll \left(\frac{X_iL_j}{|v|X_j}\right)^{\omega_1}
(X_jL_j)^{\omega_1+1}L_i
$$$$
\ll X_{j+1}^{\omega_1} X_j
L_j^{2\omega_1+1}L_i\stackrel{\eqref{eq24}}\ll
X_{j+1}^{\omega_1-(2\omega_1+1)\lambda} X_{i+1}^{-\lambda}X_j.
$$
or
$$
X_j\gg
\frac{X_iX_{i+1}^\lambda}{X_{j+1}^{1+\omega_1-(2\omega_1+1)\lambda}}.
$$

Now we substitute this estimate for $X_j$ into~\eqref{eq41} to get
$$
X_{j+1}^{1-\lambda}X_{i+1}^{1-\lambda}\gg X_i^\lambda
X_{i+1}^{\lambda^2}X_{j+1}^{(-1-\omega_1+(2\omega_1+1)\lambda)\lambda}X_i
$$$$
\Longrightarrow X_{i+1}^{1-\lambda-\lambda^2} \gg X_i^{1+\lambda}
X_{j+1}^{-1+((2\omega_1+1)\lambda-\omega_1)\lambda}.
$$
One can easily check that for $0<\lambda<\frac12$, the degree of
$X_{j+1}$ in the last expression is negative. Hence we can use the
inequality $X_{j+1}<X_i$ and Lemma~\ref{lem1} to get
$$
X_{i+1}^{1-\lambda-\lambda^2} \gg X_i^{((2\omega_1+1)\lambda
-\omega_1+1)\lambda}\gg X_{i+1}^{\frac{((2\omega_1+1)\lambda
-\omega_1+1)\lambda^2}{1-\lambda}}.
$$
That implies
$$
(1-\lambda-\lambda^2)(1-\lambda)\ge ((2\omega_1+1)\lambda
-\omega_1+1)\lambda^2.
$$

By rearranging terms, one can easily check that the last inequality
is equivalent to~\eqref{th3_eq} with $\lambda$ and $\omega_1$ in
place of $\widehat{\lambda}_3(\xi)$ and $\omega_1(\xi)$
respectively. Since $\omega_1$ can be taken arbitrarily close to
$\omega_1(\xi)$ and $\lambda$ can be taken arbitrarily close to
$\widehat{\lambda}_3(\xi)$ we finish the proof of Theorem~\ref{th3}.
\endproof

\bigskip
\noindent Dzmitry Badziahin\\ \noindent The University of Sydney\\
\noindent Camperdown 2006, NSW (Australia)\\
\noindent {\tt dzmitry.badziahin@sydney.edu.au}

\end{document}